\documentclass[11pt,a4paper,english]{amsart}
\usepackage{babel}
\usepackage[latin1]{inputenc}
\usepackage{amsmath}
\usepackage{amssymb}

\theoremstyle{plain}
\newtheorem{Thm}{Theorem}[section]

\newtheorem{Lemma}[Thm]{Lemma}

\newcommand{\la}{\lambda}

\newcommand{\dibar}{\overline {d_{i,0}}}
\newcommand{\dij}{d_{i,j}}

\theoremstyle{definition}

\newtheorem{Def}[Thm]{Definition}
\newtheorem{Def*}{Definition}

\begin{document}

\title [A proof of Parisi's conjecture]
{A proof of Parisi's conjecture on the random assignment problem}

\author{Svante Linusson}
\address{Svante Linusson\\Department of Mathematics\\
   Link{\"o}pings universitet \\
   SE-581 83 Link{\"o}ping, Sweden}
\email{linusson@mai.liu.se}

\author{Johan W{\"a}stlund}
\address{Johan W{\"a}stlund\\ Department of Mathematics \\
   Link{\"o}pings universitet \\
   SE-581 83 Link{\"o}ping, Sweden}
\email{jowas@mai.liu.se}

\date{\today}

\begin{abstract}
An assignment problem is the optimization problem of finding, in
an $m$ by $n$ matrix of nonnegative real numbers, $k$ entries, no
two in the same row or column, such that their sum is minimal.
Such an optimization problem is called a random assignment problem
if the matrix entries are random variables. We give a formula for
the expected value of the optimal $k$-assignment in a matrix where
some of the entries are zero, and all other entries are
independent exponentially distributed random variables with mean
1. Thereby we prove the formula $1+1/4+1/9+\dots+1/k^2$
conjectured by G.~Parisi for the case $k=m=n$, and the generalized
conjecture of D.~Coppersmith and G.~B.~Sorkin for arbitrary $k$,
$m$ and $n$.
\end{abstract}
\maketitle

\section{Introduction}\label{S:Intro}
The problem of minimizing the sum of $k$ elements in a matrix of
nonnegative real numbers under the condition that no two of them
may be in the same row or column is called an \emph{assignment
problem}. A set of matrix positions no two in the same row or
column is called an \emph{independent} set. An independent set of
$k$ matrix positions will also be called a $k$-assignment.

A \emph{random assignment problem}, or RAP for short, is given by
a number $k$, and an $m$ by $n$ matrix ($\min(m,n)\geq k$) of
random variables. If $P$ is a random assignment problem, we denote
by $E(P)$ the expected value of the minimal sum of an independent
set of $k$ matrix elements.

In this article we prove the following.

\begin{Thm} [Parisi's Conjecture \cite{P98}]\label{T:Parisi}
Let $P$ be the RAP where $k=m=n$ and the matrix entries are
independent exponential random variables with intensity 1.  Then
\[ E(P)=1+\frac{1}{4}+\frac{1}{9}+\cdots+\frac{1}{k^2}.
\]
\end{Thm}

We also prove the following two generalizations.

\begin{Thm} [Conjectured by D.~Coppersmith and G.~B.~Sorkin
\cite{CS98}]\label{T:CS} Let $P$ be an RAP where the matrix entries
are independent exponential random variables with intensity 1.
Then \begin{equation} \label{CSeq} E(P)
=\sum_{i+j<k}{\frac{1}{(m-i)(n-j)}}.
\end{equation}
\end{Thm}

\begin{Thm} [Conjectured in \cite{LW00}] \label{T:LW}
Let $P$ be an RAP where some matrix entries are zero and all the
other entries are independent exponential random variables with
intensity 1.  Then
\begin{equation}\label{cover}
E(P)=\frac{1}{mn}\sum_{i,j}{\frac{d_{i,j}(P)}
{\binom{m-1}{i}\binom{n-1}{j}}},
\end{equation}
where $d_{i,j}(P)$ is an integer coefficient defined in terms of
the combinatorics of the set of zeros, see Section
\ref{SS:covers}.
\end{Thm}

The identity \eqref{cover} will be referred to as the \emph{cover
formula}. When $P$ has no zeros,
$$\dij(P)=\binom{m}{i}\binom{n}{j}.$$ Hence Theorem \ref{T:LW}
$\Rightarrow$ Theorem \ref{T:CS}. To see the implication Theorem
\ref{T:CS} $\Rightarrow$ Theorem \ref{T:Parisi} (proved in
\cite{CS98}), note that if we put $m=n=k$ in \eqref{CSeq}, then
the terms for which $\gcd(k-i,k-j)=d$ sum to $1/d^2$ for
$d=1,\dots,k$.

As a consequence of Theorem \ref{T:Parisi} we also obtain a new
and completely different proof of the following theorem,
conjectured by M.~M\'ezard and G.~Parisi \cite{MP85}.

\begin{Thm} [D.~Aldous \cite{A, A00}]
\label{T:Aldous} Let $P$ be the RAP where $k=m=n$ and the matrix
entries are independent exponential random variables with intensity 1.  Then
\[ \lim_{k\to\infty}{E(P)}= \frac{\pi^2}{6}.
\]
\end{Thm}

In Section \ref{S:asymptotics} we mention some other corollaries
of Theorem \ref{T:LW}.

\subsection{Outline of the proof}

A key result, proved in Section \ref{S:row}, is a formula for the
probability that a nonzero row in an RAP is used in the optimal
$k$-assignment, see Theorem \ref{T:row}. From \cite{LW00} we know
that the probability that a nonzero element in an RAP $P$ is used
in the optimal $k$-assignment can be written $E(P)-E(P')$, where
$P'$ is obtained from $P$ by setting the matrix element in
question to zero, see Theorem \ref{T:e-d}. Therefore, the formula
for the probability that a row (or column) is used gives linear recursions for the values of the RAP's. Provided $m$ or $n$ is sufficiently large compared to $k$, this
system of linear recursions has a unique solution given by the
cover formula in Theorem \ref{T:LW}, see Section \ref{S:cover}. Finally
in Section \ref{S:rational} we prove that for fixed $k$, fixed
$m$, and a fixed set of zeros, $E(P)$ is given by a rational
function in $n$, which must then agree with the cover formula.

\subsection{Background}
Random assignment problems have attracted the attention of
researchers from physics, optimization, and probability. There are
experimental results in  \cite{Olin,PR93}. Constructive upper and
lower bounds on $E(P)$ have been given in \cite{W79, Olin, CS98,
L93, K87, GK93}. M\'ezard and Parisi \cite{MP85} used the
non-rigorous \emph{replica method} and arrived at the conjectured
limit $\pi^2/6$. This limit was subsequently established
rigorously by Aldous \cite{A, A00} using the \emph{weak
convergence} method on a weighted infinite tree model. In
this paper we continue the \emph{exact formulas}-approach inspired
by \cite{P98} and developed further in
\cite{AS02,BCR02,CS98,CS02,LW00,EES01}.

There are also interesting results on similar problems such as
finding a minimal spanning tree in a graph with random edge
weights \cite{BFM98,FM89,F85, EES01}. An intriguing question is
why the Riemann $\zeta$-function appears in the limit of both the
spanning tree and bipartite matching problems.

One tool that
we believe can be useful to a wider range of problems is
Theorem \ref{T:e-d} below and its generalization Theorem 7.3 of
\cite{LW00}.

\section{Preliminaries} \label{S:prel}
\subsection {Probabilistic preliminaries}
We say that a random variable $X$ is exponentially distributed
with intensity $a$ if $Pr(X>t)=e^{-at}$ for $t\geq 0$. The
intensity of $X$ is denoted $I(X)$. We have $E(X)=1/I(X)$.

A random assignment problem $P$ is called \emph{standard} if the
matrix entries are either zero or independent exponentially
distributed with intensity 1. A standard RAP is determined by the
numbers $k$, $m$ and $n$, and the set $Z$ of zero elements.

The following is a well-known lemma.

\begin{Lemma} Let $X_1, \dots X_n$ be independent exponential variables
with intensities $a_1,\dots,a_n$ respectively. Then the
probability that $X_i$ is minimal among $X_1, \dots X_n$ is
$$\frac{a_i}{a_1+\dots+a_n}.$$ The minimum is an exponential
variable $Y$ of intensity $a_1+\dots+a_n$ which is independent of
which variable is minimal. Under the condition that $X_i$ is
minimal, $X_i=Y$, and for $j\neq i$, we can write $X_j=Y+X_j'$,
where $X'_j$ is exponential of intensity $a_i$, and the variables
$Y$ and $X_j$ for $j\neq i$ are all independent.
\end{Lemma}

We say that a deterministic assignment problem is \emph{generic}
if no sum of a set of nonzero matrix elements is equal to the sum
of a different set of nonzero elements. In the RAP's that we
consider, the distributions of the nonzero matrix elements are
continuous. Hence all RAP's considered here will be generic with
probability 1. In the generic case, a nonzero element is used
either in every optimal $k$-assignment, or in none. Hence we can
without ambiguity speak of the probability that a certain nonzero
element is used in \emph{the} optimal $k$-assignment, without
specifying whether we take this to mean some optimal
$k$-assignment, or every optimal $k$-assignment.

Notice that this is not the case for zero elements. Even in a
generic case, there may be several different optimal
$k$-assignments, that differ in the choice of zero elements.

The following theorem, Theorem 2.10 of \cite{LW00}, is essential
for the recursion equations in Section \ref{S:cover}.

\begin{Thm} \label{T:e-d}
Suppose $P$ is a standard RAP. Suppose that the entry in position
$(i,j)$ is not zero. Let $P'$ be the standard RAP where we have
replaced the entry in position $(i,j)$ in $P$ with a zero. Then
the probability that $(i,j)$ belongs to an optimal $k$-assignment
in $P$ is $E(P)-E(P')$.
\end{Thm}

\subsection{Covers}\label{SS:covers}
We will consider sets of rows and columns in the $m\times
n$-matrix. A set $\la$ of rows and columns is said to {\em cover}
a set of zeros $Z$ if every matrix position in $Z$ is either in a
row or in a column that belongs to $\la$. A cover with $s$ rows
and columns will be called an $s$-cover. $k-1$-covers will be of
particular importance. By a {\em partial $k-1$-cover} of $Z$, we
mean a set of rows and columns which is a subset of a $k-1$-cover
of $Z$.

By a cover of an RAP we mean a cover of its set of zeros. The
\emph{cover coefficient} $d_{i,j}(P)$ is the number of partial
$k-1$-covers of $P$ with $i$ rows and $j$ columns. For this to be
nonzero, $i$ and $j$ have to be nonnegative integers with $i+j<k$.
It is convenient to regard the cover coefficient as well-defined,
but zero, for integers $i,j$ outside this range.

We say that a set $\la$ of rows and columns is an \emph{optimal
cover} of $P$, if $\la$ covers $P$, and $\la$ has minimal
cardinality among all covers of $P$. The following lemma is
well-known. For a general introduction to matching theory we refer
to \cite{LP86}.

\begin{Lemma}[Lattice property of optimal covers]
The set of optimal covers of a set of matrix positions forms a
lattice, where one of the lattice operations is taking union of
row sets and intersection of column sets, and the opposite lattice
operation is taking intersection of row sets and union of column
sets.

In particular, there is a row-maximal optimal cover containing
every row that belongs to some optimal cover, and similarly a
column-maximal optimal cover containing every column that belongs
to some optimal cover.
\end{Lemma}

An RAP can be reformulated in a setting of bipartite graphs with
random weights on the edges. An assignment is then a matching. In
this setting a cover of rows and columns is a vertex cover of the
subgraph of edges with weight zero. We will borrow terminology
from matching theory and speak of paths, etc.
referring to the corresponding graph concepts.

The following two theorems were together with Theorem \ref{T:e-d}
our main tools for computing $E(P)$ recursively in \cite{LW00}.

\begin{Thm} \label{T:mainrec} Let $A$ be a deterministic assignment problem.
If there is no set of $k$ independent zeros, then each row and
column which belongs to an optimal cover of the zeros must be used
in every optimal $k$-assignment.

Let $c$ be an optimal cover, and let $x$ be a positive real number
smaller than or equal to the minimum of the elements not covered
by $c$. Let $A'$ be obtained from $A$ by subtracting $x$ from the
elements not covered by $c$, and adding $x$ to the doubly covered
elements. Then $x$ has been subtracted from the optimal
$k$-assignment $k-|c|$ times, and consequently $$E(A) =
(k-|c|)x+E(A').$$
\end{Thm}

The previous theorem will be used for random assignment problems,
by conditioning on where the minimal non-covered element is.

\begin{Thm} \label{T:deleterow}
Suppose $P$ is a standard RAP. Suppose also that a column $c$ is
used in every $k-1$-cover of $P$. In particular, this is the case
if at least $k$ elements in $c$ are zero. Let $P'$ be the RAP with
column $c$ deleted from $P$ and $k$ and $n$ decreased by one. Then
\begin{equation}\label{eq:deleterow}
E(P)=E(P').
\end{equation}
Similarly, if a row belongs to every $k-1$-cover of $P$, it can be
deleted.
\end{Thm}

\section{Combinatorics of two optimal $k$-assignments}\label{S:manyoptimal}

In this section we prove some results that are needed in Section
\ref{S:row}. Consider a deterministic assignment problem. If there
are two different $k$-assignments $\mu$ and $\nu$, their symmetric
difference $\mu\triangle\nu:=(\mu\backslash \nu) \cup
(\nu\backslash \mu)$ will be of special importance. By a
\emph{$(\mu,\nu)$-alternating path} we mean a sequence
$(i_1,j_1),(i_2,j_2),\dots,(i_r,j_r)$ of matrix positions, where
positions belonging to $\mu$ alternate with positions belonging to
$\nu$ and $i_x=i_{x+1}$ or $j_{x}=j_{x+1}$, for all
$x=1,\dots,r-1$. The parameter $r$ will be called the
\emph{length} of the path. It is easy to see that the positions in
$\mu\triangle\nu$ will form $(\mu,\nu)$-alternating paths in the
matrix.

If two $k$-assignments $\mu$ and $\nu$ differ only at zero
positions we say that $\mu$ and $\nu$ are \emph{equivalent},
$\mu\equiv \nu$.

A deterministic assignment problem is \emph{semi-generic} if there
is exactly one nontrivial sum of nonzero elements that equals a
sum of a distinct set of nonzero elements.

\begin{Lemma} \label{L:twocomp}
Let $A$ be a deterministic assignment problem. Let $\mu$ be an
optimal $k$-assignment, and let $a$ be an element which does not
belong to $\mu$. Suppose that there is another optimal
$k$-assignment $\nu$ that contains $a$. Then there is an optimal
$k$-assignment $\nu'$ containing $a$, such that $\mu\triangle\nu'$
consists either of one path, or of two paths of odd
length.
\end{Lemma}

\begin{proof}
Let $T$ be the path of $\mu\triangle\nu$ that contains $a$.
If $T$ is of even length, then  we let $\nu' = \mu\triangle T$. If
$T$ has odd length, then it contains one more element of one of
$\mu$ and $\nu$ than of the other. Since $\mu$ and $\nu$ have the
same size, there must be another path $T'$ of
$\mu\triangle\nu$ that balances, so that $T\cup T'$ has equally
many elements from $\mu$ and $\nu$. Then we let
$\nu'=\mu\triangle(T\cup T')$. Clearly $\nu'$ has the desired
properties.
\end{proof}

\begin{Lemma} \label{L:twocomp:semigeneric}
Let $A$ be a semi-generic matrix with two inequivalent optimal
$k$-assignments $\mu$ and $\nu$. Then $\mu\triangle \nu$ has at
most two paths with some nonzero element, and if there are
two, they both have odd length.
\end{Lemma}

\begin{proof}
By Lemma \ref{L:twocomp} there is an optimal $k$-assignment $\nu'$
that contains a nonzero element not in $\mu$, such that
$\mu\triangle \nu'$ has either one path, or two paths of
odd length. Since there are only two equivalence classes of optimal
$k$-assignments, we must have $\nu'\equiv \nu$.
\end{proof}

\begin{Lemma}\label{L:tozero}
Let $A$ be a deterministic assignment problem. Let $r_1$ and $r_2$
be rows, and let $c_1$ and $c_2$ be columns. Suppose that $A(r_1,
c_1)=0$, $A(r_2, c_1)>0$, and $A(r_2, c_2)>0$. Suppose further
that there is an optimal $k$-assignment that uses $(r_2, c_2)$,
but not row $r_1$. Then there is no optimal $k$-assignment that
uses $(r_2, c_1)$.
\end{Lemma}

\begin{proof}
Let $\mu$ be an optimal $k$-assignment that uses $(r_2, c_2)$, and
suppose that there is an optimal $k$-assignment $\nu$ that uses
$(r_2, c_1)$. By Lemma \ref{L:twocomp}, we may assume that
$\mu\triangle\nu$ consists of at most two paths, and that if
there are two paths, both contain an odd number of matrix
elements.

There must be a zero element $(s_1, c_1)\in\mu$, otherwise $(r_2,
c_2)$ could be replaced by $(r_1, c_1)$ in $\mu$. Let $S$ be the
path of $(\mu\triangle\nu)\setminus\{(r_2,c_1)\}$ that
contains $(s_1, c_1)$. In other words, $S$ is the set of matrix
positions in the $(\mu,\nu)$-alternating path that starts at
$(s_1,c_1)$ and continues in the direction opposite to that of
$(r_2, c_1)$.

There must be an element $(r_1, d)$ of $\nu$ since otherwise
$(r_2,c_1)$ could be replaced by $(r_1,c_1)$ in $\nu$. Since $\mu$
does not use row $r_1$, the element $(r_1,d)$ is in a path of
$\mu\triangle\nu$ which is a path with one end belonging to $\nu$.
If the other end belongs to $\mu$, the path has an even
number of elements, and must therefore be the same as the
path containing $(r_2,c_1)$. Otherwise the other end too
belongs to $\nu$. Then the other path has both ends in $\mu$.
In either case, $S$ must end with an element of $\mu$ in a row
which is not used by $\nu$.

If the sum of the matrix entries in $\mu\cap S$ is greater than
the sum of the matrix entries in $\nu\cap S$, then the cost of the
$k$-assignment $(\mu\triangle S)\cup \{(r_1, c_1)\}$ is smaller
than the cost of $\mu$, contradicting the optimality of $\mu$.
Otherwise the sum of the matrix entries in $\mu\cap S$ is smaller
than or equal to the sum of the matrix entries in $\nu\cap S$.
Then the cost of $(\nu\triangle S)\setminus \{(r_2, c_1)\}$ is
smaller than the cost of $\nu$, again a contradiction.
\end{proof}

\begin{Lemma} \label{L:unique}
Let $A$ be a semi-generic deterministic assignment problem, and
suppose that there are two non-equivalent optimal $k$-assignments
$\mu$ and $\nu$. Assume further that
\begin{enumerate}
\item $\mu$ has a nonzero element in the last row and
\item $\nu$ does not use the last row.
\end{enumerate}
Then there is a unique row $s$ such that
\begin{enumerate}
\item $\nu$ has a nonzero element in row $s$
\item There is a $k$-assignment $\mu'\equiv\mu$ that does not use row $s$.
\end{enumerate}
\end{Lemma}

\begin{proof}
Existence: Choose a $k$-assignment $\mu'\equiv \mu$ which has as
many matrix positions as possible in common with $\nu$. Since the
last row is used by $\mu'$ but not by $\nu$, there has to be a row
$s$ which is used by $\nu$ but not by $\mu'$. If $\nu$ has a zero
element in row $s$, then $\mu'$ must contain a zero element in the
same column. By replacing this zero by the zero in row $s$, we
would obtain a $k$-assignment $\mu''\equiv\mu'$ which has one more
element in common with $\nu$, a contradiction. Therefore, the
element in row $s$ which belongs to $\nu$ must be nonzero.

Uniqueness: Assume on the contrary that there are two different
rows $s$ and $t$ that contain nonzero elements of $\nu$, say
$(s,c_1)$ and $(t,c_2)$, and two equivalent $k$-assignments
$\mu\equiv \mu'$ such that $\mu$ does not use row $s$ and $\mu'$
does not use row $t$. If $\mu$ would not use row $t$ either, then
$\mu\triangle \nu$ would have two distinct paths both with
one end in a nonzero element of $\nu$, in the elements $(s,c_1)$
and $(t,c_2)$. This contradicts Lemma \ref{L:twocomp:semigeneric}.

Hence there is a zero element $(s,c_3)$ of $\mu'$, and similarly a
zero element $(t,c_4)$ of $\mu$. Note that the symmetric
difference $\mu\triangle\mu'$ contains only zeros and therefore it
must consist of a number of $(\mu,\mu')$-alternating paths, all of
even length. Let $U$ be the one that contains $(t,c_4)$. Unless
$U$ ends at row $s$, $\mu\triangle U$ will be a $k$-assignment
equivalent to $\mu$ avoiding both rows $s$ and $t$. We have
already seen that this is impossible.

We may therefore assume that $\mu\triangle \mu'$ consists of a
single path from $(t,c_4)$ to $(s, c_3)$. Note that for every row
$r$ which is used in this path, there is a $k$-assignment
$\mu''\equiv \mu$ that avoids row $r$. This is obtained by
choosing the zeros from $\mu$ in the part of $U$ that goes towards
$(t, c_4)$, and choosing them from $\mu'$ in the part that goes
towards $(s, c_3)$.

Let $(m,d)$ be the position in the last row used by $\mu$ and let
$L$ be the $(\mu,\nu)$-alternating path containing $(m,d)$. First
note that if $L$ contains $(s,c_1)$, it has to end there and thus
be of even length, which means that it has passed through
$(t,c_2)$ first.

Case 1. $L$ intersects $U$. In this case, the first element of $L$
(starting from $(m,d)$) that belongs to $U$ must be an element of
$\mu$. Suppose that this element is in row $r$. Then let
$\mu''\equiv \mu$ avoid row $r$. It follows that $\mu''\triangle
\nu$ contains a path of even length starting at $(m,d)$, and
another path containing the nonzero element $(s,c_1)$,
contradicting Lemma \ref{L:twocomp:semigeneric}.

Case 2. $L$ does not intersect $U$. Let $S$ be the path in
$\mu\triangle\nu$ containing $(s,c_1)$, $(t,c_4)$ and $(t,c_2)$ in
this order.

We know that $U$ and $S$ intersect in $(t,c_4)$. Of the positions
in $U\cap S$, let $(r,c)$ be the one which is closest to $(s,c_3)$
in $U$, and let $\mu''\equiv \mu$ avoid row $r$. Then
$\mu''\triangle \nu$ will consist of a cycle containing the
nonzero position $(s,c_1)$, and two other paths, one
containing $(t, c_2)$ and one containing $(m,d)$, contradicting
Lemma \ref{L:twocomp:semigeneric}.
\end{proof}

\section{The probability that a row is used in the optimal
assignment}\label{S:row} Crucial for our proof is the following
formula for the probability that a certain row without zeros is
included in an optimal $k$-assignment. Let $q(P)$ denote this
probability.



\begin{Thm}[Row Inclusion Theorem] \label{T:row}
Let $P$ be a standard RAP, and let $r$ be a row without zeros. The
probability that some element in $r$ belongs to the optimal
$k$-assignment is
\begin{equation} \label{eq:row}
q(P)=\frac{1}{m}\sum_{i}{\frac{\dibar(P)}{\binom{m-1}{i}}},
\end{equation}
where $\dibar$ is the number of partial $k-1$-covers of $i$ rows
not containing the row $r$.
\end{Thm}

To simplify notation we will assume that the nonzero row in the
Row Inclusion Theorem is the {\emph last row} in the matrix. A consequence
of this formula is that the probability of using a certain row
without zeros in an optimal $k$-assignment does not change if
further zeros are introduced in a row which belongs to an optimal
cover. This observation turns out to be sufficient for the proof
of the formula.

\begin{Lemma}\label{T:invariance}
Let $P$ be a standard RAP, and let $r$ be a row that belongs to an
optimal cover. Let $P'$ be an RAP obtained from $P$ by inserting
another zero in row $r$. Suppose there is a row without zeros in
$P$. Then $$q(P)=q(P').$$
\end{Lemma}

\begin{proof} Suppose that $P$ is a standard RAP where the first row
belongs to an optimal cover, and that there is an element in the
first row, say $(1,1)$, which is not zero. Suppose further that
the last row contains no zeros. We want to show that if we replace
the element $(1,1)$ by zero, the probability that the last row is
used in the optimal $k$-assignment does not change.

Let $\Omega_P$ be the probability space of all assignments of
values to the random variables in $P$, that is, the space of all
real nonnegative $m$ by $n$ matrices that have zeros in the
positions where $P$ has zeros. If $A\in \Omega_P$, we let $q(A)$
be 1 if the last row is used by an optimal $k$-assignment, and 0
otherwise. We let $A_x$, for nonnegative real $x$, denote the
matrix obtained from $A$ by setting the entry in position $(1,1)$
to $x$.

We construct a measure preserving involution $\varphi$ on
$\Omega_P$ with the property that (except possibly on a subset of
probability zero) if $A\in \Omega_P$ is a matrix where the last
row changes between being used and not being used in the optimal
$k$-assignment when $A(1,1)$ is set to zero, then in $\varphi(A)$,
the change goes the other way. In other words,
\begin{equation}\label{qeq} q(A)-q(A_0) =
q(\varphi(A)_0)-q(\varphi(A)).\end{equation}

Let $A\in\Omega_P$. If $q(A)=q(A_0)$, we let $\varphi(A)=A$.
Otherwise notice that if we slide the element in position $(1,1)$
continuously down to zero, there can be at most one point at which
the location of the optimal $k$-assignment changes, and at this
point, the element in position $(1,1)$ goes from not being used to
being used.

At the point $A_f$ where the change occurs, the matrix is
semi-generic with two non-equivalent optimal $k$-assignments. Let
$\mu$ be the one that contains an element of the last row, and let
$\nu$ be the one that doesn't. By Lemma \ref{L:unique}, there is a
unique row $s$ such that $\nu$ contains a nonzero element in row
$s$, and so that there is a $k$-assignment $\mu'\equiv \mu$ that
contains no element of row $s$. Notice that $s$ cannot be the
first row, since the first row is used by every optimal
$k$-assignment.

We let $\varphi(A)$ be the matrix obtained from $A$ by swapping
the elements in row $s$ with the corresponding elements (elements
in the same column) in the last row, except in the columns where
row $s$ has zeros.

In the analysis of this mapping, we introduce some auxiliary
matrices. Let $A'_f$ be obtained from $A_f$ by setting the
elements in the last row that are in columns where row $s$ has
zeros, to zero. By Lemma \ref{L:tozero}, we can set these elements
as small as we please without changing the location of the optimal
assignments. Hence there will be no $k$-assignment of smaller cost
than $\mu$ and $\nu$ in $A'_f$. Then since neither $\mu$ nor $\nu$
uses any zero element in row $s$, these can of course be increased
without changing the optimality of $\mu$ and $\nu$. We let $A''_f$
be the matrix where the zero positions in row $s$ are changed to
the values of the corresponding elements in the last row of $A$.
If $A$ is generic, every optimal $k$-assignment in $A''_f$ is
equivalent either to $\mu$ or to $\nu$.

Now $\varphi(A)_f$ is obtained from $A''_f$ by swapping row $s$
with the last row. This means that in $\varphi(A)_f$, there are
exactly two equivalence classes of optimal $k$-assignments, one
that includes position $(1,1)$ and one that doesn't. But since the
last row has been swapped with another row, we have, for $x\neq
f$, $$q(\varphi(A)_x)=1-q(A_x),$$ which implies (\ref{qeq}).

In $\varphi(A)_f$, row $s$ has the property expressed in the
conclusion of Lemma \ref{L:unique}, but with the roles of $\mu$
and $\nu$ interchanged. It follows that $\varphi(\varphi(A))=A$,
and in particular that $\varphi$ is invertible. The mapping
$\varphi$ is piecewise linear, and on each piece, it is a
permutation of variables. Hence it is measure preserving. This
completes the proof.
\end{proof}

\begin{proof}[Proof of Theorem \ref{T:row}]
We now establish the Row Inclusion Theorem by an inductive
argument. Let $P$ be a standard RAP, and suppose the last row
contains no zero element. Notice that by K\"onig's theorem (see
\cite{LP86}), the formula holds whenever $P$ has a $k$-assignment
of only zeros. Suppose that the formula has been established for
every standard RAP with fewer nonzero elements than $P$.

Case 1: Some optimal cover of $P$ contains a row $r$ with at least
one nonzero element. Then every optimal $k$-assignment must use
row $r$. Let $P'$ be as $P$ but with another zero in row $r$. By
Theorem \ref{T:invariance}, the probability that the last row is
used is the same in $P$ as in $P'$. Hence we only have to show
that $\dibar(P)=\dibar(P')$ for $i=0,\dots,k-1$. We have to show
that if a set of rows can be extended to a $k-1$-cover of $P$,
then it is possible to use row $r$ in this $k-1$-cover. This will
follow if we can show that row $r$ belongs to the optimal cover of
the remaining zeros.

It suffices to show that if a row $s\neq r$ is deleted from the
matrix, row $r$ still belongs to the optimal cover of the
remaining zeros. If the deletion of $s$ does not decrease the
maximal number of independent zeros, this is obvious. Suppose
therefore that the deletion of $s$ decreases the number of
independent zeros. Then $s$ belongs to an optimal cover of $P$.
Hence the row-maximal optimal cover of $P$ contains both $r$ and
$s$. When row $s$ is deleted, the remaining rows and columns
including $r$ will constitute an optimal cover of the remaining
zeros.

Case 2: There is a row in $P$ with only zeros. Let $P'$ be
obtained from $P$ by deleting this row and decreasing $k$ and $m$
by 1. The probability that the last row is used by an optimal
assignment is clearly the same in $P$ as in $P'$. By induction,
the probability that the last row is used in $P'$ is (with $k=k(P),
m=m(P)$) $$\frac{1}{m-1}\sum_{i=0}^{k-2}
{\frac{\dibar(P')}{\binom{m-2}{i}}}.$$ We have $$ \dibar(P) =
\dibar(P')+\overline{d_{i-1,0}}(P'),$$ for every $i$. Hence
\begin{multline}
\frac{1}{m}\sum_i{\frac{\overline{d_{i,0}}(P)}{\binom{m-1}{i}}}=
\frac{1}{m}\sum_{i=0}^{k-1}{\frac{\dibar(P')+
\overline{d_{i-1,0}}(P')} {\binom{m-1}{i}}}\\ =
\frac{1}{m}\sum_i{\frac{\overline{d_{i,0}}(P')}{\binom{m-1}{i}}}+
\frac{1}{m}\sum_i{\frac{\overline{d_{i,0}}(P')}{\binom{m-1}{i+1}}}\\
=\sum_i{\overline{d_{i,0}}(P')\left(\frac{1}{m\binom{m-1}{i}}+
\frac{1}{m\binom{m-1}{i+1}}\right)}=
\frac{1}{m-1}\sum_i{\frac{\overline{d_{i.0}}(P')}{\binom{m-2}{i}}}.
\end{multline}

Case 3: There is a unique optimal cover $c$ consisting of only
columns. Then these columns will be used by the optimal
$k$-assignment. Therefore the number of elements not covered by
$c$ in the optimal $k$-assignment is independent of the random
variables in the matrix. We condition on the position of the
minimal element not covered by $c$. If in each case we subtract
this minimum from all elements not covered by $c$, the same
nonzero elements will be used in the optimal $k$-assignments. Here
we are using Theorem \ref{T:mainrec} in the special case of no
doubly covered elements. Let $P^t$ be the RAP obtained by
conditioning on the minimal element in $P$ not covered by $c$
being in row $t$, and subtracting this minimum from all elements
not covered by $c$. Then $P^t$ is a standard RAP, and the new zero
occurring in row $t$ means that row $t$ belongs to an optimal
cover of $P^t$. Hence in case the new zero is in the last row,
that row must be used in the optimal $k$-assignment, while if it
is not, we can find the probability that the last row is used by
induction. We let $P\backslash t$ be the RAP obtained by deleting
row $t$ from $P$ and decreasing $k$ by 1. Case 1 and 2 implies that the probability that the
last row is used is the same in $P^t$ as in $P\backslash t$. It
follows that the probability that the last row is used in the
optimal $k$-assignment in $P$ is given by
\begin{equation}\label{probcase3}
\frac{1}{m}+\frac{1}{m(m-1)}\sum_t{\sum_i{\frac{\overline{d_{i,0}}(P\backslash
t)}{\binom{m-2}{i}}}}.\end{equation} Here we have
$$\sum_t{\overline{d_{i,0}}(P\backslash
t)}=(i+1)\overline{d_{i+1,0}}(P),$$ since both sides are equal to
the sum, taken over all $t$, of the number of partial $k-1$-covers
of $P$ that use row $t$ (and not the last row). It follows that
(\ref{probcase3}) equals
\begin{multline}
\frac{1}{m}+\frac{1}{m(m-1)}\sum_i{\frac{(i+1)\overline{d_{i+1,0}}(P)}
{\binom{m-2}{i}}}\\=
\frac{1}{m}+\frac{1}{m}\sum_i{\frac{\overline{d_{i+1,0}}(P)}{\binom{m-1}{i+1
}}}
=\frac{1}{m}\sum_i{\frac{\overline{d_{i,0}}(P)}{\binom{m-1}{i}}}.
\end{multline} This completes the proof of Theorem \ref{T:row}.

\end{proof}

\section{Proof of the cover formula for large $m$ or large $n$}\label{S:cover}
The results of the last section will now enable us to prove that
the cover formula \eqref{cover} holds for standard RAP's whenever either $m$ or
$n$ is sufficiently large compared to $k$. We first prove that the
cover formula is consistent with the Row Inclusion Theorem.
\begin{Thm}
Let $P$ be a standard RAP where the last row contains no zeros,
and let $P_t$ be obtained from $P$ by setting the element in
column $t$ of the last row to zero. If the cover formula \eqref{cover} holds for
every $P_t$, then it holds for $P$.
\end{Thm}

\begin{proof}
By Theorem \ref{T:e-d} the probability that the element in column $t$ in the last row
belongs to an optimal $k$-assignment in $P$ is equal to
$E(P)-E(P_t)$. Hence $$nE(P)-\sum_t{E(P_t)} =
\frac{1}{m}\sum_i{\frac{\overline{d_{i,0}}(P)}{\binom{m-1}{i}}}.$$

Suppose that the cover formula holds for each $P_t$. Then

$$nE(P)=\frac{1}{mn}\sum_t{\sum_{i,j}{\frac{d_{i,j}(P_t)}{\binom{m-1}{i}
\binom{n-1}{j}}}}+
\frac{1}{m}\sum_i{\frac{\overline{d_{i,0}}(P)}{\binom{m-1}{i}}}.$$

In order to prove that the cover formula holds for $P$, it is
sufficient to prove that

\begin{equation}\label{toprove}
\sum_{i,j}{\frac{d_{i,j}(P)}{\binom{m-1}{i}\binom{n-1}{j}}}=
\frac{1}{n}\sum_t{\sum_{i,j}{\frac{d_{i,j}(P_t)}{\binom{m-1}{i}\binom{n-1}{j
}}}}+ \sum_i{\frac{\overline{d_{i,0}}(P)}{\binom{m-1}{i}}}.
\end{equation}

If we write $$\overline{d_{i,0}}(P) =
\sum_j\left(\frac{\overline{d_{i,j}}(P)}{\binom{n}{j}}-
\frac{\overline{d_{i,j+1}}(P)}{\binom{n}{j+1}}\right),$$ and fix
$i$ and $j$, we see that (\ref{toprove}) will follow from the
identity

$$\frac{d_{i,j}(P)}{\binom{n-1}{j}}=\frac{\overline{d_{i,j}}(P)}{\binom{n}{j
}}- \frac{\overline{d_{i,j+1}}(P)}{\binom{n}{j+1}}+
\frac{1}{n\binom{n-1}{j}}\sum_t{d_{i,j}(P_t)}.$$

Here, partial covers that contain the last row will contribute
$1/\binom{n-1}{j}$ to both sides. If we let
$\overline{d_{i,j}}(P_t)$ denote the number of partial
$k-1$-covers of $P_t$ with $i$ rows and $j$ columns that do not
contain the last row, it only remains to show that

$$\frac{\overline{d_{i,j}}(P)}{\binom{n-1}{j}}=\frac{\overline{d_{i,j}}(P)}{
\binom{n}{j}}- \frac{\overline{d_{i,j+1}}(P)}{\binom{n}{j+1}}+
\frac{1}{n\binom{n-1}{j}}\sum_t{\overline{d_{i,j}}(P_t)},$$ or
equivalently
$$j\overline{d_{i,j}}(P)+(j+1)\overline{d_{i,j+1}}(P)=
\sum_t{\overline{d_{i,j}}(P_t)}.$$

Here the first term of the left hand side counts the partial
$k-1$-covers of $P_t$ that contain column $t$, while the second
term counts those that don't. Hence (\ref{toprove}) holds.
\end{proof}

Next we show that the cover formula is consistent with removing a
column which belongs to every $k-1$-cover.

\begin{Thm}
Let $P$ be a standard RAP, and suppose that the first column
belongs to every $k-1$-cover.  Let $P'$ be the RAP obtained from
$P$ by deleting the first column and decreasing $k$ and $n$ by 1. If the
cover formula holds for $P'$, then it holds for $P$. \end{Thm}

\begin{proof}
By Theorem \ref{T:deleterow} $E(P)=E(P')$. We have (with $n=n(P)$ and $k=k(P)$)
\begin{multline}
\frac{1}{mn}\sum_{i,j}{\frac{d_{i,j}(P)}{\binom{m-1}{i}\binom{n-1}{j}}}
=\frac{1}{mn}\sum_{i,j}{\frac{d_{i,j-1}(P')+d_{i,j}(P')}{\binom{m-1}{i}
\binom{n-1}{j}}}\\
=\frac{1}{mn}\sum_{i,j}{\left(\frac{d_{i,j}(P')}{\binom{m-1}{i}\binom{n-1}{j
}}+ \frac{d_{i,j}(P')}{\binom{m-1}{i}\binom{n-1}{j+1}}\right)}=
\frac{1}{m(n-1)}\sum_{i,j}{\frac{d_{i,j}(P')}{\binom{m-1}{i}\binom{n-2}{j}}}\\
=E(P')=E(P).
\end{multline}

\end{proof}
\begin{Thm}
If $P$ is a standard RAP such that $\max(m,n)>(k-1)^2$, then
$$E(P)=\frac{1}{mn}\sum_{i,j}{\frac{d_{i,j}(P)}
{\binom{m-1}{i}\binom{n-1}{j}}}.$$
\end{Thm}

\begin{proof}
We show that if $m>(k-1)^2$, then a standard RAP which contains a
zero in each row must have a column which belongs to every
$k-1$-cover. The theorem then follows by induction on the number
of nonzero elements of $P$. It suffices to consider the case that
$P$ has exactly one zero in each row. Then the total number of
zeros is greater than $(k-1)^2$. Either there is a column with $k$
or more zeros, which must then belong to every $k-1$-cover, or
there is an independent set of $k$ zeros, and in this case there
is no $k-1$-cover.

By symmetry, the formula also holds if $n>(k-1)^2$.
\end{proof}

\section{Rationality of $E(P)$ as a function of $n$}
\label{S:rational} We know from Section \ref{S:cover} that for fixed $k$, whenever $m$ or $n$
is large, the cover formula \eqref{cover} for standard RAP's holds. In order to
prove that the formula holds for smaller values of $m$ and $n$, it
is therefore sufficient to show that if $k$, $m$ and the zero
positions are fixed, and we let $P_n$ be the standard RAP with $n$
columns, then there is a rational function in the variable $n$
that gives the value of $E(P_n)$ for every $n$ which is at least
as large as $k$ and the number of columns with zeros. This
rational function must then be equal to the one given by the cover
formula. We prove this by induction over a class of RAP's which
includes not only standard RAP's.

In an \emph{exponential} RAP, all matrix elements are linear
combinations with nonnegative rational coefficients of a set $X_1,
\dots, X_p$ of independent exponentially distributed random
variables.
If a variable $X_i$ in an exponential RAP has intensity 1, and
occurs in one and only one matrix position, and this matrix entry
is equal to $X_i$, then the variable is called a \emph{standard
variable}, and the matrix position where it occurs is called a
\emph{standard position}.

We introduce the concept of an \emph{RAP-sequence}. The idea is to
treat a set of similar RAP's with different number of columns in a
uniform way, in order to prove that there is a rational expression
in the number of columns that gives the value of each RAP in the
set.

\begin{Def}
We say that a linear function $f(x)=ax+b$ in one variable is
\emph{$k$-positive}, iff $f(x)>0$ whenever $x\geq k$, or
equivalently, if $a\geq 0$ and $b>-ak$.
\end{Def}

Obviously, the sum of two or more $k$-positive functions is
$k$-positive.

\begin{Def}
An \emph{RAP-sequence} is a sequence $P_n$ ($n\geq n_0$) of
exponential RAP's satisfying the following:
\begin{enumerate}
\item
Each $P_n$ is an exponential RAP in a set of variables
$X_{n,1},\dots, X_{n, p(n)}$.

\item The numbers $m$ and $k$ are uniform, that is $m(P_i)=m(P_j)$
and $k(P_i)=k(P_j)$ for all $i$ and $j$, while $P_n$ has $n$
columns.

\item In the first $n_0$ columns, the matrix elements of the various RAP's
differ only in that the first index of the variables is changed.
In other words, the coefficient of $X_{n,i}$ in a matrix element
in the first $n_0$ columns of $P_n$ is equal to the coefficient of
$X_{n_0, i}$ in the same position in $P_{n_0}$.

\item Beyond column $n_0$, $P_n$ has only standard elements.
\end{enumerate}

Moreover, we say that the RAP-sequence is \emph{well-behaved} if
\begin{enumerate}
\item For each $i$ such that the variables $X_{n,i}$ occur in the
first $n_0$ columns, there is a $k$-positive linear polynomial
$f_i(n)$ such that the intensity of $X_{n,i}$ is $f_i(n)$.

\item
Nonzero nonstandard elements occur only in columns which belong to
the column-maximal optimal cover of the zeros, or equivalently,
columns that intersect every maximal set of independent zeros.
\end{enumerate}
\end{Def}

\begin{Thm}[Rationality Theorem] \label{T:rationality}
Suppose $P$ is a well-behaved RAP-sequence. Then there is a
rational function $f(x)$ in one variable such that
\begin{enumerate}
\item If $x$ is a zero of the denominator of $f$, then $x<k$.
\item $E(P_n)=f(n)$ for every $n\geq n_0(P)$.
\end{enumerate}
\end{Thm}

\begin{Def}
We say that a linear combination $u_1$ of variables
$X_1,\dots,X_p$ is smaller than another if for each $X_i$, the
coefficient is smaller in $u_1$ than in $u_2$. We say that $u_1$
and $u_2$ are \emph{incomparable} if neither of them is smaller
than the other, in other words, if each of them has a higher
coefficient than the other for some variable.

We say that a linear combination is \emph{potentially minimal} in
a set of linear combinations, if it is smaller than or
incomparable with every other.
\end{Def}

Let $c$ be the row-maximal optimal cover of $P$. 
We prove the rationality theorem by induction on a number of
parameters, in the following order:
\begin {enumerate}
\item The size of the largest independent set of zeros. An
RAP-sequence is considered simpler if it has a larger set of
independent zeros.

\item The number of rows in $c$.
If the number of independent zeros are equal, the
RAP-sequence with fewer rows belonging to the row-maximal optimal
cover is simpler.

\item The set of potentially minimal nonstandard elements not
covered by $c$. If 1. and 2. are equal, an RAP-sequence is
considered simpler if fewer of the nonstandard elements not
covered by $c$ are potentially minimal.

\item If 1--3 are equal, and there are two incomparable
nonstandard elements not covered by $c$, then an RAP-sequence is
simpler if there are fewer variables with different coefficients
in the first two (in lexicographic order, say) incomparable
potentially minimal nonstandard elements not covered by $c$.

\item If 1--3 are equal, and there is a minimal non-covered
nonstandard element, then an RAP-sequence is simpler if the number
of variables occurring in this element is smaller.
\end{enumerate}

Let $P$ be a well-behaved RAP-sequence, and suppose that the
rationality theorem holds for every simpler well-behaved
RAP-sequence with the same values of $m$ and $k$. We may of course assume that
$|c|<k$. Since $P$ is well-behaved, in each row all but at most
$k-1$ elements are standard and not covered by $c$.

We show that $E(P)$ can be expressed in terms of rational
functions in $n$, and values of simpler well-behaved
RAP-sequences.

Case 1: There are two or more non-covered incomparable nonstandard
elements. Let $u_1$ and $u_2$ be the first two (in lexicographic order). 
We choose $i$ and $j$ such that the
coefficient of $X_{n,i}$ is greater in $u_1$, and the coefficient
of $X_{n,j}$ is greater in $u_2$. Let the coefficients be $a_1$,
$a_2$, $b_1$ and $b_2$ so that $u_1 = a_1X_{n,i} + b_1X_{n,j}
+\dots$ and $u_2=a_2X_{n,i} +b_2X_{n,j}+\dots$.

Let $Q$ and $R$ be the RAP-sequences obtained by conditioning on
$(a_1-a_2)X_{n,i}$ being smaller or greater than
$(b_2-b_1)X_{n,j}$, respectively.

The intensities of $(a_1-a_2)X_{n,i}$ and $(b_2-b_1)X_{n,j}$ are
$f_i(n)/(a_1-a_2)$ and $f_j(n)/(b_2-b_1)$ respectively. The
probability of $(a_1-a_2)X_{n,i}$ being smaller than
$(b_2-b_1)X_{n,j}$ is
$$\frac{f_i(n)/(a_1-a_2)}{f_i(n)/(a_1-a_2)+f_j(n)/(b_2-b_1)},$$
and similarly, the probability of $(b_2-b_1)X_{n,j}$ being smaller
than $(a_1-a_2)X_{n,i}$ is
$$\frac{f_j(n)/(b_2-b_1)}{f_i(n)/(a_1-a_2)+f_j(n)/(b_2-b_1)}.$$

Therefore, \begin{equation} \label{case1} E(P_n)=
\frac{f_i(n)E(Q_n)/(a_1-a_2)+f_j(n)E(R_n)/(b_2-b_1)}
{f_i(n)/(a_1-a_2)+f_j(n)/(b_2-b_1)}.\end{equation} We show that
$Q$ and $R$ can be regarded as well-behaved RAP-sequences. If we
condition on $(a_1-a_2)X_{n,i}$ being smaller than
$(b_2-b_1)X_{n,j}$, then we can write $$(a_1-a_2)X_{n,i}=Y_n$$ and
$$(b_2-b_1)X_{n,j}=Y_n+Z_n,$$ where $Y_n$ and $Z_n$ are
independent exponentially distributed variables. The intensities
are given by
$$I(Y_n)=\frac{f_i(n)}{a_1-a_2}+\frac{f_j(n)}{b_2-b_1}$$ and
$$I(Z_n)=\frac{f_j(n)}{b_2-b_1},$$ both of which are $k$-positive.

When replacing $X_{n,i}$ and $X_{n,j}$ by the new variables $Y_n$
and $Z_n$, only the nonstandard elements are affected. Either
$u_1$ is smaller than $u_2$ in $Q$, or at least the number of
variables with distinct coefficients is smaller than in $P$, since
$u_1$ and $u_2$ will get the same coefficient for $Y_n$.
For any two matrix elements that satisfied $u_3\le u_4$ 
before we conditioned on which one of
$X_{n,i}$ and $X_{n,j}$ is smallest, this inequality will still hold. 

Hence $Q$, and similarly $R$, are simpler than $P$. By induction,
it follows that (\ref{case1}) gives a rational expression for
$E(P)$, where the denominator is nonzero for $n\geq k$.

Case 2: There are no two non-covered incomparable nonstandard
elements. Then either there is among the non-covered nonstandard
elements a minimal one, or there are no non-covered nonstandard
elements. We can treat these slightly different cases in the same
way.

If there is a minimal non-covered nonstandard element, let
$aX_{n,i}$ be a term occurring in this matrix entry. We let $S_n$
be a set of random variables consisting of $aX_{n,i}$ and all the
non-covered standard variables (if there is no non-covered
nonstandard element, we let $S_n$ consist only of the non-covered
standard elements). There are at least $n-k+1$ non-covered
standard elements in each non-covered row. By grouping together
the standard elements in each row, we can write the total
intensity of $S_n$ in a uniform way as a sum of $k$-positive
terms. Therefore the intensity of the minimum $Y_n$ of the terms
in $S_n$ is $k$-positive as a function of $n$.

We condition on the minimal element in $S_n$. We can then replace
the terms in $S_n$ by new variables $Y_n$ and $Z_{n,i}$ where
$Y_n=\min(S_n)$ and $Z_{n,i}$ are the differences between the
remaining terms in $S_n$ and the minimum.

Since all non-covered nonstandard elements contain the variable
$X_{n,i}$ with a coefficient of at least $a$, we can by Theorem \ref{T:mainrec}
subtract the
minimum of $S_n$ from every non-covered element, and add it to the
doubly covered elements.

We then get
\begin{equation}\label{case2}
E(P_n)=\frac{k-|c|}{I(Y_n)}+\frac{1}{I(Y_n)}\sum_{t\in
S_n}{I(t)E(Q_n(t))},
\end{equation}
where $Q_n(t)$ are the new RAP's obtained by conditioning on the
term $t$ being smallest, and performing the change of variables
and subtraction of the minimum, and $I(t)/I(Y_n)$ is the
corresponding probability. In order to write this in a uniform way
for the different values of $n$, we group together the cases where
the minimum occurs in a particular row and beyond column $n_0$.

By permuting columns, we may assume that in those cases, the
minimum always occurs in column $n_0+1$. The probability for each
such case is $(n-n_0)/I(Y_n)$. In this way, (\ref{case2}) will
contain the same terms for each $n>n_0$.

If we let $S'_n$ be the subset of $S_n$ consisting of variables
occurring in the first $n_0$ columns, and we let $S''_n$ be the
set of standard variables in column $n_0+1$, we get

\begin{equation}\label{case2'}
E(P_n)=\frac{k-|c|}{I(Y_n)}+\frac{1}{I(Y_n)}\sum_{t\in
S'_n}{I(t)E(Q_n(t))}+\frac{n-n_0}{I(Y_n)}\sum_{t\in
S''_n}{I(t)E(Q_n(t))}.
\end{equation}

For $n=n_0$, the second sum will be empty. However, provided we
can show that $Q(t)$ is a well-behaved RAP-sequence simpler than
$P$, it will follow by induction that the denominator of $E(Q(t))$
does not vanish for $n\leq k$. Hence the rational expression that
occurs when multiplying the probability $(n-n_0)/I(S_n)$ with the
expression for $E(Q(t))$ will vanish for $n=n_0$. Therefore the
rational expression that results will give the correct value of
$E(P)$ also for $n=n_0$.

It remains to show that $Q(t)$ is a well-behaved RAP-sequence, and
that it is simpler than $P$. We first show that $Q(t)$ is
well-behaved. We have already seen that the intensities of the
variables occurring in $Q(t)$ are $k$-positive. We therefore turn
to the distribution of nonstandard elements. Possibly, there are
some new nonstandard elements among the doubly covered elements.
If there is no new zero among the non-covered elements, the
columns of the doubly covered positions of course belong to the
column-maximal optimal cover. Suppose that a new zero occurs among
the non-covered elements. If the new zero is in a column that
belongs to the column-maximal cover of $P$, then this is still the
column-maximal cover of the new set of zeros. If the new zero is
in a column that does not belong to the column-maximal cover, then
since it is also in a row that does not belong to the row-maximal
cover, there must be an independent set of zeros in $Q(t)$ which
is larger than the largest independent set of zeros in $P$.
Therefore, the column-maximal cover in $P$, extended with the
column where the new zero has occurred, will be an optimal cover
in $Q(t)$. Hence in any case, the columns of the column-maximal
cover in $P$ belong to the column-maximal cover of $Q(t)$.

We now show that $Q(t)$ is simpler than $P$. If $Q(t)$ contains a
new zero, then either it gives a larger independent set of zeros,
or it has to be covered by a column in the row-maximal cover in
$Q(t)$. In either case, $Q(t)$ is simpler than $P$. If no new zero
occurs, this must be because the minimal term in $S_n$ was a term
occurring in the minimal non-covered nonstandard element. In this
case, the number of variables in this element will decrease, again
making $Q(t)$ simpler than $P$.

Hence (\ref{case2}) gives a rational expression for $E(P)$ whose
denominator is non-vanishing for $n\geq k$.

This completes the proof of Theorem \ref{T:rationality}.

We are finally able to give a proof of Theorem \ref{T:LW}, which
we restate.

\noindent
{\bf Theorem \ref {T:LW}}
\emph {If $P$ is a standard RAP, then}
$$E(P)=\frac{1}{mn}\sum_{i,j}{\frac{d_{i,j}(P)}
{\binom{m-1}{i}\binom{n-1}{j}}}.$$

\begin{proof}
A standard RAP with $n_0$ columns can be extended to a
well-behaved RAP-sequence $P_n$ by inserting more columns without
zeros. Hence Theorem \ref{T:rationality} shows that there is a
rational function $f$ giving the expected value of $P_n$. Since
the cover coefficient $d_{i,j}(P_n)$ can be expressed uniformly as
a polynomial in $n$, the cover formula gives a rational function
in $n$ which takes the same values as $f$ on the infinitely many
integers $n>\max(n_0, (k-1)^2)$. Hence the cover formula must
agree with $f$, and give the value of $E(P)$.
\end{proof}

\section{Asymptotic results}\label{S:asymptotics}
Besides giving a new proof of Aldous' $\zeta(2)$-limit theorem,
Theorem \ref{T:LW} also makes the conjectured limits of
\cite{LW00} rigorous. It follows already from the conjecture of
Coppersmith and Sorkin, now Theorem \ref{T:CS}, that if $\alpha,
\beta\geq 1$, then as $k\to \infty$, the value of the optimal
$k$-assignment in a $[\alpha k]$ by $[\beta k]$-matrix of
exp(1)-variables converges to $$\int_\Delta{\frac{dx
dy}{(\alpha-x)(\beta-y)}},$$ where $\Delta$ is the triangle with
vertices in $(0,0)$, $(1,0)$ and $(0,1)$. For instance, when
$\alpha=1$ and $\beta=2$ the limit is equal to
$$\frac{\pi^2}{12}-\frac{(\log{2})^2}{2}.$$

It also follows that the value of a standard RAP with zeros in a
region which is scaled up with $k$, $m$ and $n$ will converge to a
similar integral. In particular it is shown in \cite{LW00} that
the cover formula implies that the limit value of a standard
square RAP with zeros outside an inscribed circle is equal to
$\pi^2/24$.

A result presented as a conjecture in \cite{Olin}, and proved in
\cite{A00}, states that in the case $k=m=n$ with no zeros, as
$n\to \infty$, the probability that the smallest element in a row
belongs to the optimal assignment converges to $1/2$. We can now
give an exact formula for this probability for finite $k$, $$
\frac{1}{2}+\frac{1}{2k}.$$ In the case $k=m$, the probability
that the smallest element in a particular row belongs to the
optimal assignment is equal to the probability that the smallest
element in the entire matrix does. For arbitrary $k$, $m$, and
$n$, the probability that the smallest element in the matrix
belongs to the optimal $k$-assignment is $$1-\frac{k(k-1)}{2mn}.$$

These theorems are obtained by combining the results of
\cite{LW00} and Theorem \ref{T:LW}.

\end{document}